\documentclass{amsart}

\usepackage{amssymb,amsmath,amsthm,latexsym,stmaryrd}

 \newtheorem{thm}{Theorem}[section]
 \newtheorem{cor}[thm]{Corollary}

 \newtheorem{prop}[thm]{Proposition}
 \theoremstyle{definition}
 
 \theoremstyle{remark}

\newcommand{\ie}{i.e.\ }
\newcommand{\f}{\varphi}
\newcommand{\tg}{\tilde{g}}
\newcommand{\n}{\nabla}
\newcommand{\nn}{\tilde{\n}}
\newcommand{\M}{(M,\allowbreak{}\f,\allowbreak{}\xi,\allowbreak{}\eta,\allowbreak{}g)}
\newcommand{\G}{\mathcal{G}}
\newcommand{\I}{\mathcal{I}}

\newcommand{\R}{\mathbb R}
\newcommand{\C}{\mathbb C}
\newcommand{\X}{\mathfrak X}
\newcommand{\F}{\mathcal{F}}

\newcommand{\ta}{\theta}
\newcommand{\om}{\omega}
\newcommand{\lm}{\lambda}

\newcommand{\al}{\alpha}
\newcommand{\bt}{\beta}

\newcommand{\dsfrac}{\displaystyle\frac}

\newcommand{\Span}{\operatorname{Span}}
\newcommand{\Ker}{\operatorname{Ker}}
\newcommand{\Id}{\operatorname{Id}}
\newcommand{\ad}{\operatorname{Ad}}

\newcommand{\thmref}[1]{Theorem~\ref{#1}}

\newcommand{\corref}[1]{Corollary~\ref{#1}}

\begin{document}
%
%
%
%
%
%
%
%
%
\title{Lie Groups as 3-Dimensional Almost Contact B-Metric Manifolds}

\author[H. Manev, D. Mekerov]{Hristo Manev, Dimitar Mekerov
}

\address[H. Manev]{Medical University of Plovdiv, Faculty of Pharmacy, 
Department of Pharmaceutical Sciences,   15-A Vasil Aprilov
Blvd.,   Plovdiv 4002,   Bulgaria;}  \address{Paisii
Hilendarski University of Plovdiv,   Faculty of Mathematics and
Informatics,   Department of Algebra and Geometry,   236
Bulgaria Blvd.,   Plovdiv 4027,   Bulgaria}
\email{hmanev@uni-plovdiv.bg}
\address[D. Mekerov]{Paisii Hilendarski University of Plovdiv, 
Faculty of Mathematics and Informatics,   Department of Algebra
and Geometry,   236 Bulgaria Blvd.,   Plovdiv 4027, 
Bulgaria} \email{mircho@uni-plovdiv.bg}





\begin{abstract}
The object of investigations are almost contact B-metric
structures on 3-dimensional Lie groups considered as smooth
manifolds. There are established the existence and some geometric
characteristics of these manifolds in all basic classes. An
example is given as a support of obtained results.
\end{abstract}

\maketitle

\section*{Introduction}

The differential geometry of the almost contact metric manifolds
is well studied (e.g. \cite{Blair}). The beginning of the
investigations on the almost contact manifolds with B-metric is
given by \cite{GaMiGr}. These manifolds are the odd-dimensional
counterpart of the almost complex manifolds with Norden metric
\cite{GaBo,GrMeDj}, where the almost complex structure acts as an
anti-isometry on the metric. Almost contact B-metric manifolds are
investigated and studied in
\cite{GaMiGr,Man4,Man31,ManGri1,ManGri2,ManIv13,ManIv14,NakGri2}.

An object of special interest is the case of the lowest dimension
of the considered manifolds. Curvature properties of the
3-dimensional manifolds of this type are studied in
\cite{Man33,ManNak15,NakMan16}.

The basic problem considered in this work is the existence and the
geometric characteristics of almost contact B-metric structures on
3-dimensional Lie algebras. The main result  is a construction of
the considered structures and their characterization in all basic
classes of the almost contact B-metric manifolds in the
classification given in \cite{GaMiGr}.

The present paper is organized as follows. In
Sect.~\ref{sect-prel} we recall
some facts about the almost contact B-metric manifolds. %
%
In Sect.~\ref{sect-Lie-mfds} we construct different types of
contact structures on Lie algebras  in dimension 3.
In Sect.~\ref{sect-Exms} we consider an example in relation with
the investigations in the previous section.


\section{Preliminaries}\label{sect-prel}

\subsection{Almost contact manifolds with
B-metric}\label{sect-mfds}

Let $(M,\f,\xi,\eta,g)$ be an almost contact manifold with
B-met\-ric or an \emph{almost contact B-metric manifold}, \ie $M$
is a $(2n+1)$-dimensional differentiable manifold with an almost
contact structure $(\f,\xi,\eta)$ consisting of an endomorphism
$\f$ of the tangent bundle, a vector field $\xi$, its dual 1-form
$\eta$ as well as $M$ is equipped with a pseu\-do-Rie\-mannian
metric $g$  of signature $(n+1,n)$, such that the following
algebraic relations are satisfied \cite{GaMiGr}:
\begin{equation}\label{str}
\begin{array}{c}
\f\xi = 0,\qquad \f^2 = -\Id + \eta \otimes \xi,\qquad
\eta\circ\f=0,\qquad \eta(\xi)=1,\\[4pt]
g(\f x, \f y) = - g(x,y) + \eta(x)\eta(y),
\end{array}
\end{equation}
where $\Id$ is the identity and $x$, $y$ are elements of the
algebra $\X(M)$ of the smooth vector fields on $M$.

Further $x$, $y$, $z$ will stand for arbitrary elements of $\X(M)$
or vectors in the tangent space $T_pM$ of $M$ at an arbitrary
point $p$ in $M$.

The associated metric $\tg$ of $g$ on $M$ is defined by
\(\tg(x,y)=g(x,\f y)+\eta(x)\eta(y)\).  The manifold
$(M,\f,\xi,\eta,\tg)$ is also an almost contact B-metric manifold.
The signature of both metrics $g$ and $\tg$ is necessarily
$(n+1,n)$. By $\n$ and $\nn$ we denote the Levi-Civita connection
of $g$ and $\tg$, respectively.

A classification of almost contact B-metric manifolds is given in
\cite{GaMiGr}. This classification, consisting of eleven basic
classes $\F_1$, $\F_2$, $\dots$, $\F_{11}$,  is made with respect
to the tensor  $F$ of type (0,3) defined by
\begin{equation*}\label{F=nfi}
F(x,y,z)=g\bigl( \left( \nabla_x \f \right)y,z\bigr).
\end{equation*}
The following properties are valid:
\begin{equation*}\label{F-prop}
F(x,y,z)=F(x,z,y)=F(x,\f y,\f z)+\eta(y)F(x,\xi,z)
+\eta(z)F(x,y,\xi).
\end{equation*}

The intersection of the basic classes is the special class $\F_0$
determined by the condition $F(x,y,z)=0$. Hence $\F_0$ is the
class of almost contact B-metric manifolds with $\n$-parallel
structures, \ie $\n\f=\n\xi=\n\eta=\n g=\n \tg=0$.

If $\left\{e_i;\xi\right\}$ $(i=1,2,\dots,2n)$ is a basis of
$T_pM$ and $\left(g^{ij}\right)$ is the inverse matrix of
$\left(g_{ij}\right)$, then with $F$ are associated 1-forms
$\theta$, $\theta^*$, $\omega$, called \emph{Lee forms}, defined
by:
\begin{equation*}\label{t}
\theta(z)=g^{ij}F(e_i,e_j,z),\quad \theta^*(z)=g^{ij}F(e_i,\f
e_j,z), \quad \omega(z)=F(\xi,\xi,z).
\end{equation*}

Further, we use the following characteristic conditions of the
basic classes \cite{Man8}:
\begin{equation}\label{Fi}
\begin{array}{rl}
\F_{1}: &F(x,y,z)=\frac{1}{2n}\bigl\{g(x,\f y)\ta(\f z)+g(\f x,\f
y)\ta(\f^2 z)\\[4pt]
&\phantom{F(x,y,z)=\frac{1}{2n}\bigl\{}+g(x,\f z)\ta(\f y)+g(\f
x,\f z)\ta(\f^2 y)
\bigr\};\\[4pt]
\F_{2}: &F(\xi,y,z)=F(x,\xi,z)=0,\quad\\[4pt]
              &F(x,y,\f z)+F(y,z,\f x)+F(z,x,\f y)=0,\quad \ta=0;\\[4pt]
\F_{3}: &F(\xi,y,z)=F(x,\xi,z)=0,\quad\\[4pt]
              &F(x,y,z)+F(y,z,x)+F(z,x,y)=0;\\[4pt]
\F_{4}: &F(x,y,z)=-\frac{1}{2n}\ta(\xi)\bigl\{g(\f x,\f y)\eta(z)+g(\f x,\f z)\eta(y)\bigr\};\\[4pt]
\F_{5}: &F(x,y,z)=-\frac{1}{2n}\ta^*(\xi)\bigl\{g( x,\f y)\eta(z)+g(x,\f z)\eta(y)\bigr\};\\[4pt]
\F_{6}: &F(x,y,z)=F(x,y,\xi)\eta(z)+F(x,z,\xi)\eta(y),\quad \\[4pt]
                &F(x,y,\xi)= F(y,x,\xi)=-F(\f x,\f y,\xi),\quad \ta=\ta^*=0; \\[4pt]
\F_{7}: &F(x,y,z)=F(x,y,\xi)\eta(z)+F(x,z,\xi)\eta(y),\quad \\[4pt]
                &F(x,y,\xi)=-F(y,x,\xi)=-F(\f x,\f y,\xi); \\[4pt]
\F_{8}: &F(x,y,z)=F(x,y,\xi)\eta(z)+F(x,z,\xi)\eta(y),\\[4pt]
                &F(x,y,\xi)=F(y,x,\xi)=F(\f x,\f y,\xi); \\[4pt]
\F_{9}: &F(x,y,z)=F(x,y,\xi)\eta(z)+F(x,z,\xi)\eta(y),\\[4pt]
                &F(x,y,\xi)=-F(y,x,\xi)=F(\f x,\f y,\xi); \\[4pt]
\F_{10}: &F(x,y,z)=F(\xi,\f y,\f z)\eta(x); \\[4pt]
\F_{11}:
&F(x,y,z)=\eta(x)\left\{\eta(y)\om(z)+\eta(z)\om(y)\right\}.
\end{array}
\end{equation}

Let now consider the case of the lowest dimension of the manifold
$\M$, \ie $\dim{M}=3$.

Let us denote the components $F_{ijk}=F(e_i,e_j,e_k)$ of the
tensor $F$ with respect to a \emph{$\f$-basis}
$\left\{e_i\right\}_{i=0}^2=\left\{e_0=\xi,e_1=e,e_2=\f e\right\}$
of $T_pM$, which satisfies the following conditions
\begin{equation}\label{gij}
g(e_0,e_0)=g(e_1,e_1)=-g(e_2,e_2)=1,\quad g(e_i,e_j)=0,\; i\neq j
\in \{0,1,2\}.
\end{equation}

According to \cite{HM}, the components of the Lee forms with
respect to a $\f$-basis with conditions \eqref{gij} are
\begin{equation*}\label{t3}
\begin{array}{lll}
\ta_0=F_{110}-F_{220},\qquad & \ta_1=F_{111}-F_{221},\qquad &\ta_2=F_{112}-F_{211},\\[4pt]
\ta^*_0=F_{120}+F_{210},\qquad &\ta^*_1=F_{112}+F_{211},\qquad &\ta^*_2=F_{111}+F_{221},\\[4pt]
\om_0=0,\qquad &\om_1=F_{001},\qquad &\om_2=F_{002}.
\end{array}
\end{equation*}

Then, if $F_s$ $(s=1,2,\dots,11)$ are the components of $F$ in the
corresponding basic classes $\F_s$ and $x=x^ie_i$, $y=y^je_j$,
$z=z^ke_k$ are arbitrary vectors we have \cite{HM}:

\begin{equation}\label{Fi3}
\begin{array}{rl}
&F_{1}(x,y,z)=\left(x^1\ta_1-x^2\ta_2\right)\left(y^1z^1+y^2z^2\right),\\[4pt]
&\qquad\ta_1=F_{111}=F_{122},\qquad \ta_2=-F_{211}=-F_{222}; \\[4pt]
&F_{2}(x,y,z)=F_{3}(x,y,z)=0;
\\[4pt]
&F_{4}(x,y,z)=\frac{1}{2}\ta_0\Bigl\{x^1\left(y^0z^1+y^1z^0\right)
-x^2\left(y^0z^2+y^2z^0\right)\bigr\},\\[4pt]
&\qquad \frac{1}{2}\ta_0=F_{101}=F_{110}=-F_{202}=-F_{220};\\[4pt]
&F_{5}(x,y,z)=\frac{1}{2}\ta^*_0\bigl\{x^1\left(y^0z^2+y^2z^0\right)
+x^2\left(y^0z^1+y^1z^0\right)\bigr\},\\[4pt]
&\qquad \frac{1}{2}\ta^*_0=F_{102}=F_{120}=F_{201}=F_{210};\\[4pt]
&F_{6}(x,y,z)=F_{7}(x,y,z)=0;\\[4pt]
&F_{8}(x,y,z)=\lm\bigl\{x^1\left(y^0z^1+y^1z^0\right)
+x^2\left(y^0z^2+y^2z^0\right)\bigr\},\\[4pt]
&\qquad \lm=F_{101}=F_{110}=F_{202}=F_{220};\\[4pt]
&F_{9}(x,y,z)=\mu\bigl\{x^1\left(y^0z^2+y^2z^0\right)
-x^2\left(y^0z^1+y^1z^0\right)\bigr\},\\[4pt]
&\qquad \mu=F_{102}=F_{120}=-F_{201}=-F_{210};\\[4pt]
&F_{10}(x,y,z)=\nu x^0\left(y^1z^1+y^2z^2\right),\qquad
\nu=F_{011}=F_{022};\\[4pt]
&F_{11}(x,y,z)=x^0\bigl\{\left(y^1z^0+y^0z^1\right)\om_{1}
+\left(y^2z^0+y^0z^2\right)\om_{2}\bigr\},\\[4pt]
&\qquad \om_1=F_{010}=F_{001},\qquad \om_2=F_{020}=F_{002}.
\end{array}
\end{equation}

Obviously, the class of 3-dimensional almost contact B-metric
manifolds is
\begin{equation}\label{thm-3D}
\F_1 \oplus \F_4 \oplus \F_5 \oplus \F_8 \oplus \F_9 \oplus
\F_{10} \oplus \F_{11}.
\end{equation}

\subsection{Curvature properties of almost contact B-metric manifolds}\label{sec-curv}

 Let $R=\left[\n,\n\right]-\n_{[\ ,\ ]}$ be the
curvature (1,3)-tensor of $\nabla$ and the corresponding curvature
$(0,4)$-tensor be denoted by the same letter: $R(x,y,z,w)$
$=g(R(x,y)z,w)$. The following properties are valid:
\begin{equation}\label{R}
\begin{array}{c}
    R(x,y,z,w)=-R(y,x,z,w)=-R(x,y,w,z), \\[4pt]
R(x,y,z,w)+R(y,z,x,w)+R(z,x,y,w)=0.
\end{array}
\end{equation}

On any $\F_0$-mani\-fold we have
\begin{equation}\label{Kel}
R(x,y,\f z,\f w)=-R(x,y,z,w).
\end{equation}

Any  (0,4)-tensor having the properties \eqref{R} is called a
\emph{curvature-like} tensor. If such a tensor has the K\"ahler
property \eqref{Kel} it is called a \emph{K\"ahler} tensor.

From \cite{ManNak15} it is known that every K\"ahler tensor on a
3-dimensional almost contact B-metric manifold is zero and
therefore every 3-dimen\-sion\-al $\F_0$-manifold is flat, \ie $R
= 0$.

 The Ricci
tensor $\rho$ and the scalar curvature $\tau$ for $R$ as well as
their associated quantities are defined respectively by
\begin{equation*}
\begin{array}{ll}
    \rho(y,z)=g^{ij}R(e_i,y,z,e_j),\qquad &
    \tau=g^{ij}\rho(e_i,e_j),\\[4pt]
    \rho^*(y,z)=g^{ij}R(e_i,y,z,\f e_j),\qquad &
    \tau^*=g^{ij}\rho^*(e_i,e_j).
\end{array}
\end{equation*}

Further, we use the notation $g\owedge  h$ for the Kulkarni-Nomizu
product of two (0,2)-tensors, \ie
\[
\begin{split}
\left(g\owedge h\right)(x,y,z,w)&=g(x,z)h(y,w)-g(y,z)h(x,w)\\[4pt]
&+g(y,w)h(x,z)-g(x,w)h(y,z).
\end{split}
\]
Obviously, $g\owedge h$ has the properties of $R$ in \eqref{R}
when $g$ and $h$ are symmetric.

%
%
%

It is well known that the curvature tensor  has the following form
on any 3-dimensional manifold:
\begin{equation}\label{R3}
R 
=-g\owedge \left(\rho-\frac{\tau}{4} g\right).
\end{equation}

Every non-degenerate 2-plane $\beta$ with respect to $g$ in
$T_pM$, $p \in M$, has the following sectional curvature
\begin{equation*}\label{sect}
k(\beta;p)=-\frac{2R(x,y,y,x)}{g\owedge g},
\end{equation*}
where $\{x,y\}$ is a basis of $\beta$.

A 2-plane $\beta$ is said to be a \emph{$\f$-holomorphic section}
(respectively, a \emph{$\xi$-section}) if $\beta= \f\beta$
(respectively, $\xi \in \beta$).

\section{Lie Groups as 3-Dimensional Almost Contact B-Metric
Manifolds}\label{sect-Lie-mfds}

Let $L$ be a 3-dimensional real connected Lie group, and
$\mathfrak{l}$ be its Lie algebra with a basis
$\{E_{0},E_{1},E_{2}\}$ of left invariant vector fields. Let an
almost contact structure $(\f,\xi,\eta)$ be defined by
\begin{equation}\label{strL}
\begin{array}{l}
\f E_0=0,\quad \f E_1=E_{2},\quad \f E_{2}=- E_1,\quad \xi=
E_0,\quad \\[4pt]
\eta(E_0)=1,\quad \eta(E_1)=\eta(E_{2})=0
\end{array}
\end{equation}
and let $g$ be a pseudo-Riemannian metric such that
\begin{equation}\label{gL}
\begin{array}{l}
  g(E_0,E_0)=g(E_1,E_1)=-g(E_{2},E_{2})=1, \\[4pt]
  g(E_0,E_1)=g(E_0,E_2)=g(E_1,E_2)=0.
\end{array}
\end{equation}
Then we have a 3-dimensional almost contact B-metric manifold
$(L,\f,\xi,\eta,g)$.

The corresponding Lie algebra $\mathfrak{l}$
 is determined as follows:
\begin{equation}\label{lie}
\left[E_{i},E_{j}\right]=C_{ij}^k E_{k}, \quad i, j, k \in \{0,1,2\}.
\end{equation}

Applying the Jacobi identity
\[
\bigl[[E_i,E_j],E_k\bigr]+\bigl[[E_j,E_k],E_i\bigr]+\bigl[[E_k,E_i],E_j\bigr]=0,
\]
from the nine commutation coefficients $C_{ij}^k$ remain six which
could be chosen as parameters, so the three coefficients with
different indices are expressed by the six parameters (if the
denominators are non-zero) as follows:
\begin{equation*}\label{Jac}
\begin{array}{c}
C_{12}^0=\dsfrac{C_{12}^1 C_{01}^0+C_{12}^2 C_{02}^0}{C_{01}^1 +
C_{02}^2},\qquad C_{02}^1=\dsfrac{C_{01}^1 C_{02}^0-C_{12}^1
C_{02}^2}{C_{01}^0-C_{12}^2},
\\[4pt]
C_{01}^2=\dsfrac{C_{01}^0 C_{02}^2+C_{01}^1 C_{12}^2}{C_{02}^0 +
C_{12}^1}.
\end{array}
\end{equation*}

The Koszul equality
\begin{equation}\label{Kosz}
2g\left(\n_{E_i}E_j,E_k\right)
=g\left([E_i,E_j],E_k\right)+g\left([E_k,E_i],E_j\right)
+g\left([E_k,E_j],E_i\right)
\end{equation}
implies the following formula for the components
$F_{ijk}=F(E_i,E_j,E_k)$, $i,j,\allowbreak{}k\in\{0,1,2\}$, of the
tensor $F$:
\begin{equation*}\label{Fijk}
\begin{split}
2F_{ijk}=g\left([E_i,\f E_j]-\f [E_i,E_j],E_k\right)
&+g\left(\f [E_k,E_i]-[\f E_k,E_i],E_j\right)\\[4pt]
&+g\left([E_k,\f E_j]-[\f E_k,E_j],E_i\right).
\end{split}
\end{equation*}

Then we have
\begin{equation}\label{FijkC}
\begin{array}{l}
\begin{array}{ll}
F_{111}=F_{122}=2C_{12}^1,\quad &
F_{211}=F_{222}=2C_{12}^2,\\[4pt]
F_{120}=F_{102}=-C_{01}^1,\quad &
F_{020}=F_{002}=-C_{01}^0,\\[4pt]
F_{210}=F_{201}=-C_{02}^2,\quad &
F_{010}=F_{001}=C_{02}^0,\\[4pt]
\end{array}\ \\[4pt]
\begin{array}{l}
F_{110}=F_{101}=\frac12 \left(C_{12}^0-C_{01}^2+C_{02}^1\right),\\[4pt]
F_{220}=F_{202}=\frac12 \left(C_{12}^0+C_{01}^2-C_{02}^1\right),\\[4pt]
F_{011}=F_{022}=C_{12}^0+C_{01}^2+C_{02}^1.
\end{array}
\end{array}
\end{equation}
For the Lee forms we obtain
\begin{equation}\label{titiC}
\begin{array}{l}
\begin{array}{l}
\ta_{0}=-C_{01}^2+C_{02}^1,\\[4pt]
\ta_{1}=2C_{12}^1,\\[4pt]
\ta_{2}=-2C_{12}^2,\\[4pt]
\end{array}\quad
\begin{array}{l}
\ta^*_{0}=-C_{01}^1-C_{02}^2,\\[4pt]
\ta^*_{1}=2C_{12}^2,\\[4pt]
\ta^*_{2}=2C_{12}^1,\\[4pt]
\end{array}\quad
\begin{array}{l}
\om_{0}=0,\\[4pt]
\om_{1}=C_{02}^0,\\[4pt]
\om_{2}=-C_{01}^0.
\end{array}
\end{array}
\end{equation}

Using \eqref{Fi3}, \eqref{FijkC} and \eqref{titiC}, we deduce the following
\begin{thm}\label{thm-Fi-L}
The manifold $(L,\f,\xi,\eta,g)$ belongs to the basic class $\F_s$
($s \in \{1,4,5,8,9,10,11\}$) if and only
if the corresponding Lie algebra $\mathfrak{l}$ is determined by
the following commutators:
\begin{equation}\label{Fi-L}
\begin{array}{ll}
\F_1:\; &[E_0,E_1]=[E_0,E_2]=0, \quad [E_1,E_2]=\al E_1+\bt E_2;
\\[4pt]
\F_4:\; &[E_0,E_1]=\al E_2, \quad [E_0,E_2]=-\al E_1, \quad
[E_1,E_2]=0;
\\[4pt]
\F_5:\; &[E_0,E_1]=\al E_1, \quad [E_0,E_2]=\al E_2, \quad
[E_1,E_2]=0;
\\[4pt]
\F_8:\; &[E_0,E_1]=\al E_2, \quad [E_0,E_2]=\al E_1, \quad
[E_1,E_2]=-2\al E_0;
\\[4pt]
\F_9:\; &[E_0,E_1]=\al E_1, \quad [E_0,E_2]=-\al E_2, \quad
[E_1,E_2]=0;
\\[4pt]
\F_{10}:\; &[E_0,E_1]=\al E_2, \quad [E_0,E_2]=\al E_1, \quad
[E_1,E_2]=0;
\\[4pt]
\F_{11}:\; &[E_0,E_1]=\al E_0, \quad [E_0,E_2]=\bt E_0, \quad
[E_1,E_2]=0,
\end{array}
\end{equation}
where $\al$, $\bt$ are arbitrary real parameters. Moreover, the
relations of $\al$ and $\bt$ with the non-zero components
$F_{ijk}$ in the different basic classes $\F_s$ from \eqref{Fi3} are the
following:
\begin{equation*}\label{Fi-L-alpha}
\begin{array}{ll}
\F_1:\quad \al=\frac12 \ta_1,\quad \bt=\frac12 \ta_2; \qquad
&\F_4:\quad \al=\frac12 \ta_0;\\[4pt]
\F_5:\quad \al=-\frac12 \ta^*_0; \qquad
&\F_8:\quad \al=-\lm; \\[4pt]
\F_9:\quad \al=-\mu; \qquad
&\F_{10}:\quad \al=\frac12 \nu;   \\[4pt]
\F_{11}:\quad \al=-\om_2, \quad \bt=\om_1.\qquad &
\end{array}
\end{equation*}
\end{thm}

Let us remark that an $\F_{0}$-manifold is obtained when the Lie
algebra is Abelian, \ie all commutators are zero. Further, we omit
this special case.

\subsection{Structures of special type on the considered manifolds}


The metric $g$ is called \emph{Killing} (or,
\emph{$\ad$-invariant}) if
\begin{equation}\label{g-Kil}
g([x,y],z)=g(x,[y,z]).
\end{equation}

Bearing in mind \eqref{gL}, \eqref{lie} and \eqref{g-Kil}, we
obtain that $g$ is Killing if and only if the equalities
$C_{12}^0=-C_{01}^2=-C_{02}^1$ are valid and the rest of the
commutation coefficients are zero. Then, according to
\eqref{FijkC} and \eqref{titiC}, we get the following non-zero
components
$-F_{011}=-F_{022}=2F_{110}=2F_{101}=2F_{220}=2F_{202}=C_{12}^0$
and $\ta=\ta^*=\om=0$. Therefore, using \eqref{Fi3}, we establish
that $F=F_8+F_{10}$ and $2\lm=-\nu$.
Thus, it is valid the following
\begin{thm}
The metric $g$ is Killing if and only if $(L,\f,\xi,\eta,g)$
belongs to the subclass of $\F_8\oplus\F_{10}$ determined by the
condition $2\lm=-\nu$.
\end{thm}

In a similar way we obtain that $\tg$ is Killing, i.e.
$\tg([x,y],z)=\tg(x,[y,z])$, if and only if
$C_{12}^0=-C_{01}^1=C_{02}^2$ and the others $C_{ij}^k$ are zero.
Then from \eqref{FijkC} we have zero Lee forms and
$F_{120}=F_{102}=-F_{210}=-F_{201}=2F_{110}=2F_{101}
=2F_{220}=2F_{202}=F_{011}=F_{022}=C_{12}^0$, \ie $2\lm=\mu=\nu$,
bearing in mind \eqref{Fi3}. Therefore, we get
\begin{thm}
The metric $\tg$ is Killing if and only if $(L,\f,\xi,\eta,g)$
belongs to the subclass of $\F_8\oplus\F_9\oplus\F_{10}$
determined by the condition $2\lm=\mu=\nu$.
\end{thm}

The structure $\f$ is called \emph{bi-invariant}, if $\f[x, y] =
[x, \f y]$. This condition is satisfied if and only if
$F_{011}=F_{022}=2F_{110}=2F_{101}=2F_{220}=2F_{202}=C_{12}^0$,
\ie $2\lm=\nu$.
  Therefore, it is valid the following
\begin{thm}
The structure $\f$ is bi-invariant if and only if
$(L,\f,\xi,\eta,g)$ belongs to the subclass of $\F_8\oplus\F_{10}$
determined by the condition $2\lm=\nu$.
\end{thm}

The structure $\f$ is called \emph{Abelian}, if $[\f x, \f y] =
[x, y]$. This condition is valid if and only if
$F_{011}=F_{022}=2F_{110}=2F_{101}=2F_{220}=2F_{202}=C_{12}^0$,
\ie $2\lm=\nu$, and $F_{111}=F_{122}=2C_{12}^1=\ta_1$,
$F_{211}=F_{222}=2C_{12}^2=-\ta_2$. Therefore, we deduce
\begin{thm}
The structure $\f$ is Abelian if and only if $(L,\f,\xi,\eta,g)$
belongs to the subclass of $\F_1\oplus\F_8\oplus\F_{10}$
determined by the condition $2\lm=\nu$.
\end{thm}

The vector field $\xi$ is called \emph{Killing}, if the Lie
derivative with respect to $\xi$ of the metric $g$ is zero, \ie
$\mathcal{L}_{\xi}g=0$. Then the manifold $(L,\f,\xi,\eta,g)$ is
called a \emph{K-contact B-metric manifold}. This condition is
equivalent to $\left(\n_x\eta\right)y+\left(\n_y\eta\right)x=0$,
which is satisfied for any manifold in
$\F_1\oplus\F_2\oplus\F_3\oplus\F_7\oplus\F_{8}\oplus\F_{10}$ for
arbitrary dimension, according to \eqref{Fi}. Then, bearing in
mind \eqref{thm-3D} -- the class of the considered manifolds of
dimension 3, we get
\begin{thm}
The vector field $\xi$ is Killing if and only if
$(L,\f,\xi,\eta,g)\in\F_1\oplus\F_{8}\oplus\F_{10}$, \ie the class
of the 3-dimensional K-contact B-metric manifolds is
$\F_1\oplus\F_{8}\oplus\F_{10}$.
\end{thm}

\subsection{Curvature properties of the considered manifolds}

Using \eqref{Kosz} and \eqref{Fi-L}, we obtain the covariant
derivatives of $E_i$ with respect to $\n$. The non-zero ones of
them are the following:
\[
\begin{array}{ll}
\F_1: &\n_{E_1}E_1=\al E_2,\;  \n_{E_1}E_2=\al E_1, \;
\n_{E_2}E_1=-\bt E_2, \; \n_{E_2}E_2=-\bt E_1;
\\[4pt]
\F_4: &\n_{E_1}E_0=-\al E_2,\;  \n_{E_2}E_0=\al E_1;
\\[4pt]
\F_5: &\n_{E_1}E_1=-\n_{E_2}E_2=\al E_0,\;  \n_{E_1}E_0=-\al
E_1,\; \n_{E_2}E_0=-\al E_2;
\\[4pt]
\F_8: &\n_{E_1}E_2=-\n_{E_2}E_1=-\al E_0,\;  \n_{E_1}E_0=-\al
E_2,\; \n_{E_2}E_0=-\al E_1;
\\[4pt]
\F_9: &\n_{E_1}E_1=\n_{E_2}E_2=\al E_0,\;  \n_{E_1}E_0=-\al E_1,\;
\n_{E_2}E_0=\al E_2;
\\[4pt]
\F_{10}: &\n_{E_0}E_1=\al E_2,\; \n_{E_0}E_2=\al E_1;
\\[4pt]
\F_{11}: &\n_{E_0}E_0=-\al E_1+\bt E_2,\; \n_{E_0}E_1=\al E_0,
\n_{E_0}E_2=\bt E_0.
\end{array}
\]

By virtue of the latter equalities and bearing in mind the
definition  equalities of the curvature tensor $R$, the Ricci
tensor $\rho$, the $*$-Ricci tensor $\rho^*$, the scalar curvature
$\tau$, its $*$-scalar curvature $\tau^*$ and the sectional
curvature $k_{ij}$ of the non-degenerate 2-plane $\{E_i,E_j\}$
\begin{equation}\label{Rrhotauk}
\begin{array}{c}
R_{ijkl}=g\left(\n_{E_i}\n_{E_j}E_k-\n_{E_j}\n_{E_i}E_k
-\n_{[E_i,E_j]}E_k,E_l\right),
\\[4pt]
\rho_{jk}=R_{0jk0}+R_{1jk1}-R_{2jk2}, \qquad
\rho^*_{jk}=R_{1jk2}+R_{2jk1},
\\[4pt]
\tau=\rho_{00}+\rho_{11}-\rho_{22},\qquad \tau^*=2\rho_{12},
\qquad k_{ij}=-2R_{ijij}/(g\owedge g)_{ijij},
\end{array}
\end{equation}
from \thmref{thm-Fi-L} we obtain the following

\begin{thm}\label{thm-res}
The $\F_{10}$-manifolds $(L,\f,\xi,\eta,g)$ are flat and the non-zero components of $R$,
$\rho$, $\rho^*$ and the non-zero values of $\tau$, $\tau^*$,
$k_{ij}$ for the manifolds $(L,\f,\xi,\eta,g)$ from the rest of the basic classes are the following:
\begin{equation}\label{res}
\begin{array}{ll}
\F_1:\; &R_{1212}=\rho_{11}=-\rho_{22}=\rho^*_{12}=\rho^*_{21}=\frac12 \tau=k_{12} 
=
\al^2-\bt^2;\;  \\[4pt]
\F_4:\;
&-R_{0101}=R_{0202}=\frac12
\rho_{00}=\rho_{11}=-\rho_{22} \\[4pt]
&\phantom{-R_{0101}}=\frac14 \tau
=k_{01}=k_{02}=\al^2;\;\\[4pt]
\F_5:\;
&-R_{0101}=R_{0202}=R_{1212}=\frac12\rho_{00}=\frac12\rho_{11}
=-\frac12\rho_{22} \\[4pt]
&\phantom{-R_{0101}}=\rho^*_{12}=\rho^*_{21} =\frac16 \tau=k_{12}=k_{01}=k_{02}=-\al^2;\\[4pt]
\F_8:\; &R_{0101}=-R_{0202}=R_{1212}=-\frac12\rho_{00}=\rho^*_{12}=\rho^*_{21}  \\[4pt]
&\phantom{R_{0101}}=-\frac12\tau=k_{12}=-k_{01}=-k_{02}=\al^2;\;  \\[4pt]
\F_9:\; &R_{0101}=-R_{0202}=R_{1212}=-\frac12\rho_{00}=\rho^*_{12}=\rho^*_{21} \\[4pt]
&\phantom{R_{0101}}=-\frac12\tau=k_{12}=-k_{01}=-k_{02}=\al^2;\\[4pt]
%
%
\F_{11}:\; &R_{0101}=-\rho_{11}=-k_{01}=\al^2, \quad  R_{0202}=-\rho_{22}=k_{02}=\bt^2, \\[4pt]
&R_{0120}=\rho_{12}=\frac12 \tau^*=-\al\bt, \quad
\rho_{00}=\frac12 \tau=-\al^2+\bt^2.
\end{array}
\end{equation}
\end{thm}

Let us remark, on the 3-dimensional almost contact B-metric
manifold with the basis $\{E_{0},E_{1},E_{2}\}$ defined by
\eqref{strL}, we have two basic $\xi$-sections $\{E_{0},E_{1}\}$,
$\{E_{0},E_{2}\}$ and one basic $\f$-holomorphic section
$\{E_{1},E_{2}\}$.

According to \thmref{thm-res}, we obtain
\begin{thm}\label{thm-char}
For the $\F_{s}$-manifolds $(L,\f,\xi,\eta,g)$, which do not belong to $\F_{0}$, the
following propositions are valid:
\begin{enumerate}
    \item Only the $\F_{1}$-manifolds have a
    K\"ahler tensor $R$;
      \item An $\F_{1}$-manifold is flat if and only if
    $\al^2=\bt^2$;
\item All $\F_{8}$-manifolds and $\F_{9}$-manifolds
    have curvature tensors of the same form;
\item The curvature tensor of any $\F_{s}$-manifold ($s=4,11$)
    has the property $R(x,y,\f z,\f w)=0$;
\item Every $\F_{s}$-manifold ($s=4,11$) has a vanishing $*$-Ricci tensor;
    \item Every $\F_{s}$-manifold ($s=4,8,9$) has a positive scalar
    curvature;
  \item An $\F_{1}$-manifold (resp., $\F_{11}$-manifold) has a positive scalar
    curvature if and only if $\al^2>\bt^2$ (resp., $\al^2<\bt^2$);
\item Every $\F_{5}$-manifold has a negative scalar curvature;
\item  An $\F_{1}$-manifold
    (resp., $\F_{11}$-manifold) has a negative scalar
    curvature if and only if $\al^2<\bt^2$ (resp., $\al^2>\bt^2$);
   \item Every $\F_{s}$-manifold ($s=1,4,5,8,9$) has a vanishing $*$-scalar
    curvature;
\item  An $\F_{11}$-manifold has a vanishing $*$-scalar
    curvature if and only if $\al\bt=0$ for $\al\neq\bt$;
\item      An $\F_{11}$-manifold
    has a positive (resp., negative) $*$-scalar
    curvature if and only if $\al\bt<0$ (resp., $\al\bt>0$);
\item Every $\F_{1}$-manifold  has vanishing sectional curvatures of the $\xi$-sect\-ions;
   \item  Every $\F_{4}$-manifold  has positive sectional curvatures of
    the basic $\xi$-sect\-ions;
\item  Every $\F_{s}$-manifold ($s=5,8,9$) has negative sectional curvatures
    of the basic $\xi$-sections;
\item An $\F_{1}$-manifold has a vanishing scalar
    curvature of the basic $\f$-holo\-mor\-phic section if and only if the manifold is
    flat;
\item  Every $\F_{s}$-manifold ($s=4,11$) has a vanishing scalar
    curvature of the basic $\f$-holomorphic section;
\item Every $\F_{s}$-manifold
    ($s=8,9$) (resp., $\F_{5}$-manifold) has a positive (resp., negative)
    scalar curvature of the basic $\f$-holomorphic section;
\item An $\F_{1}$-manifold has a positive (resp., negative) scalar
    curvature of the basic $\f$-holomorphic section
    if and only if $\al^2>\bt^2$ (resp., $\al^2<\bt^2$).
\end{enumerate}
\end{thm}

Taking into account \eqref{res}, we get the following
\begin{cor}\label{cor-Fi3-rho}
The form of the Ricci tensor on  $(L,\f,\xi,\eta,g)$ in the respective basic class is:
\begin{equation*}\label{Fi3-rho}
\begin{array}{llll}
\F_1: & \rho=\frac{\tau}{2}\left(g-\eta\otimes\eta\right); \qquad
& \F_4: & \rho=\frac{\tau}{4}\left(g+\eta\otimes\eta\right);
\\[4pt]
\F_5: & \rho=\frac{\tau}{3}g; \quad & \F_8: &
\rho=\tau(\eta\otimes\eta);
\\[4pt]
\F_9: & \rho=\tau(\eta\otimes\eta); \quad & \F_{11}: &
\rho=\rho(\f\cdot,\f\cdot)+\frac{\tau}{2}g-\tau^*g^*,
\end{array}
\end{equation*}
where $g^*(\cdot,\cdot)=g(\cdot,\f \cdot)= \tg - \eta\otimes\eta$.
\end{cor}

Using \corref{cor-Fi3-rho} and  \eqref{R3}, we obtain for the
following
\begin{cor}\label{cor-Fi3-R}
The form of the curvature tensor on  $(L,\f,\xi,\eta,g)$ in the respective basic
class is:
\begin{equation*}\label{Fi3-R} %
\begin{array}{llll}
\F_1: & R=-\frac{\tau}{4}\left(g\owedge
g\right)+\frac{\tau}{2}\left(g\owedge
\left(\eta\otimes\eta\right)\right);
\\[4pt]
\F_4: & R=-\frac{\tau}{4}\left(g\owedge
\left(\eta\otimes\eta\right)\right);
\\[4pt]
\F_5: & R=-\frac{\tau}{12}\left(g^*\owedge g^*\right)
-\frac{\tau}{6} \left( g\owedge
(\eta\otimes\eta)\right); \\[4pt]
\F_8: & R=\frac{\tau}{4}\left(g^*\owedge g^*\right)
-\frac{\tau}{2} \left(g\owedge (\eta\otimes\eta)\right);
\\[4pt]
\F_9: & R=\frac{\tau}{4}\left(g^*\owedge g^*\right)
-\frac{\tau}{2} \left(g\owedge (\eta\otimes\eta)\right); \\[4pt]
\F_{11}: & R=-\rho \owedge (\eta\otimes\eta).
\end{array}
\end{equation*}
\end{cor}

Let  $(V,\f,\xi,\eta,g)$ be a vector space endowed with an almost contact structure $(\f,\xi,\eta)$.
The decomposition $x=-\f^2x+\eta(x)\xi$ generates the projectors
$h$ and $v$ on $(V,\f,\xi,\eta)$ determined by $h(x)=-\f^2x$ and
$v(x)=\eta(x)\xi$ and having the properties $h\circ h =h$, $v\circ
v=v$, $h\circ v=\allowbreak{}v\circ h=0$. Therefore, we have the
orthogonal decomposition $V=h(V)\oplus v(V)$. Then we obtain the
corresponding linear operators $\ell_1,\ell_2, \ell_3$ in the space $V^*\times V^*$ of the
tensors $S$ of type (0,2) over $(V,\f,\xi,\eta)$ defined in
\cite{AlGa}:
\begin{equation*}\label{ell}
\begin{array}{c}
\ell_1(S)(x,y)=S(h(x),h(y)),\quad \ell_2(S)(x,y)=S(v(x), v(y)),\\[4pt]
\ell_3(S)(x,y)=S(v(x), y)+S(x, v(y))-2S(v(x), v(y)).
\end{array}
\end{equation*}

The structure group $\G$ of $\M$ is
determined by $\G=\mathcal{O}(n;\mathbb{C})\times\I$, where $\I$
is the identity on $\Span(\xi)$ and
$\mathcal{O}(n;\mathbb{C})=\mathcal{GL}(n;\mathbb{C})\cap
\mathcal{O}(n,n)$.
%

It is known  \cite{Man-diss}, the following orthogonal
decomposition which is invariant with respect to the action of the
structure group  $\mathcal{G}$:
\begin{equation*}
V^*\times V^*=\ell_1(V^*\times V^*)\oplus \ell_2(V^*\times V^*)
\oplus \ell_3(V^*\times V^*),
\end{equation*}
where
\begin{equation*}
\ell_i(V^*\times V^*)=\left\{S\in (V^*\times V^*)\ |\
S=\ell_i(S)\right\},\quad i=1,2,3.
\end{equation*}

Since the metrics $g$ and $\tg$ belong to $V^*\times V^*$, then
their components in the three orthogonal subspaces are:
\begin{equation*}\label{ell-g}
\begin{array}{lll}
\ell_1(g)=-g(\f\cdot,\f\cdot)=g-\eta\otimes\eta, \quad
&\ell_2(g)=\eta\otimes\eta,
\quad &\ell_3(g)=0,\\[4pt]
\ell_1(\tg)=g(\cdot,\f\cdot)=\tg-\eta\otimes\eta, \quad
&\ell_2(\tg)=\eta\otimes\eta, \quad &\ell_3(\tg)=0.
\end{array}
\end{equation*}

On any almost contact B-metric manifold there are a pair of
B-metrics $(g,\tg)$ which restrictions on $H=\Ker(\eta)$ coincide
with their $\ell_1$-images and are the real part and the imaginary
part of the complex Riemannian metric $g^{\C}$, \ie
$g^{\mathbb{C}}=g|_H+i\tg|_H=\ell_1(g)+i\ell_1(\tg)$. Moreover,
the restrictions of $g$ and $\tg$ on $H^{\bot}=\Span(\xi)$
coincide with their $\ell_2$-images and $\eta\otimes\eta$, \ie
$g|_{H^{\bot}}=\tg|_{H^{\bot}}=\ell_2(g)=\ell_2(\tg)=\eta\otimes\eta$.

It is well known that a manifold is called \emph{Einstein} if the
Ricci tensor is proportional to the metric tensor, \ie $\rho=\lm
g$, $\lm\in\R$.

Then, bearing in mind the above arguments, it is reasonable to
call a manifold $M$ whose Ricci tensor satisfies the condition
\begin{equation}\label{e}
\rho=\lm g + \mu \tg + \nu \eta\otimes\eta
\end{equation}
 an
\emph{$\eta$-complex-Einstein manifold}. The case of the
$\eta$-complex-Einstein Sasaki-like manifolds is considered in
\cite{IMM}. If $\mu=0$ we call $M$ an \emph{$\eta$-Einstein
manifold}. On the almost contact metric manifolds in \cite{BoGaMa}
is studied the $\eta$-Einstein geometry as a class of
distinguished Riemannian metrics on the contact metric manifolds.
It is well-known \cite{YaKo} that $\lm$ and $\nu$ are constants
if the almost contact metric manifold is K-contact and has
dimension greater than 3.


It is meaningful to consider the Einstein property in the separate
components $\ell_i(\rho)$ and $\ell_i(\rho^*)$, $i=1,2,3$.

In \cite{ManNak15} the following notions related to the Einstein
condition are introduced. An almost contact B-metric manifold $\M$
is called \emph{contact-Einstein} if the Ricci tensor  has the
form $\rho=\lm g|_H + \mu \tg|_H + \nu \eta\otimes\eta$, where
$\lm,\mu,\nu\in\R$. A contact-Einstein manifold is called
\emph{$h$-Einstein} (resp., \emph{$v$-Einstein})  if $\rho=\lm
g|_H + \mu \tg|_H$ (resp., $\rho=\nu \eta\otimes\eta$). An
$h$-Einstein manifold is called \emph{$\f$-Einstein} (resp.,
\emph{$*$-Einstein}) if $\rho=\lm g|_H$ (resp., $\rho= \mu
\tg|_H$).

According to \corref{cor-Fi3-rho}, we obtain
\begin{prop}
The manifold $\M$ is:
\begin{enumerate}
    \item $\f$-Einstein if it belongs to $\F_1$;
    \item $\eta$-Einstein if it belongs to $\F_4$;
    \item Einstein if it belongs to $\F_5$;
    \item $v$-Einstein if it belongs to $\F_8\oplus\F_9$.
\end{enumerate}
\end{prop}

\section{An Example}
\label{sect-Exms}


In \cite{HM} is introduced a real connected Lie group $L$ as a
manifold from the class $\F_9\oplus\F_{10}$. We consider now this
example for $\dim L=3$ and give geometrical characteristics
in relation with some previous results in this work.
In this case $L$ is a 3-dimensional real connected Lie group and
its associated Lie algebra with a global basis $\{E_{0},E_{1},
E_{2}\}$ of left invariant vector fields on $L$ is defined by
\begin{equation*}\label{com}
    [E_0,E_1]=-a_1E_1-a_{2}E_{2},\quad
    [E_0,E_{2}]=-a_{2}E_1+a_{1}E_{2},\quad
[E_1,E_2]=0,
\end{equation*}
where $a_1, a_{2}$ are real constants. Actually, $L$ is a family
of Lie groups depending on the pair of real parameters.
An invariant almost contact structure is defined by
\eqref{strL} and
$g$ is a pseudo-Riemannian metric such that \eqref{gL} are
valid. Thus, because of \eqref{str}, the induced 3-dimensional
manifold $(L,\f, \xi, \eta, g)$ is an almost contact B-metric
manifold.

Let us remark that in \cite{Ol} the same Lie group with the same
almost contact structure but equipped with a compatible Riemannian
metric is studied as an almost cosymplectic manifold.

The non-zero components $F_{ijk}$ of the tensor $F$ are following:
\[
F_{011}=F_{022}=-2a_2,\quad F_{102}=F_{120}=-F_{201}=-F_{210}=a_1.
\]
Thus, it follows that $(L,\f, \xi, \eta, g)\in \F_9\oplus\F_{10}$,
bearing in mind \eqref{Fi3} for $\mu=a_1$, $\nu=-2a_2$.

Further, we continue with examination of this family of almost
contact B-metric manifolds.

Using \eqref{Kosz}, we obtain the covariant derivatives with
respect to $\n$ of $E_i$ for the manifold $(L,\f, \xi, \eta, g)$:
\[
\begin{array}{c}
\n_{E_0}E_1=\frac{\nu}{2} E_2,\quad  \n_{E_0}E_2=\frac{\nu}{2}
E_1,\quad \n_{E_1}E_0=\mu E_1,\quad   \n_{E_2}E_0=-\mu E_2,
\\[4pt]
\n_{E_1}E_1=\n_{E_2}E_2=-\mu E_0.
\end{array}
\]

By virtue of the latter equalities and bearing in mind the
definition  equalities \eqref{Rrhotauk}, we get:
\begin{equation}\label{Rrhotauk-Ex}
\begin{array}{c}
R_{0110}=-R_{0220}=R_{1221}=-\mu^2, \qquad R_{0120}=\mu\nu,
\\[4pt]
\rho_{00}=-2\rho_{22}=-2\mu^2, \qquad
\rho^*_{00}=2\rho_{12}=2\mu\nu,
\\[4pt]
\tau=-2\mu^2,\qquad \tau^*=2\mu\nu,
\\[4pt]
k_{01}=k_{02}=-\mu^2,\qquad k_{12}=\mu^2
\end{array}
\end{equation}

We obtain the form of the Ricci tensor $\rho$:
\begin{equation}\label{rho-Ex}
\rho=\left(\tau+\frac{\tau^*}{2}\right)\eta\otimes\eta-\frac{\tau^*}{2}\tg.
\end{equation}
Consequently, bearing in mind \eqref{e}, $(L,\f, \xi, \eta, g)$ is
an $\eta$-complex-Einstein manifold with $\lm=0$.

According to \eqref{R3} and \eqref{rho-Ex}, we obtain the
following form of the curvature tensor:
\begin{equation}\label{R-Ex}
R=\frac{1}{4}g\owedge \Bigl(\tau
g+2\tau^*\tg-2\left(2\tau+\tau^*\right)\eta\otimes\eta\Bigr).
\end{equation}

If $\mu=0$, $\nu\neq0$, \ie the manifold belongs to $\F_{10}$,
then $\tau=\tau^*=0$. Therefore,  \eqref{R-Ex} implies $R=0$.

If $\nu=0$, $\mu\neq0$, \ie the manifold belongs to $\F_9$, then
$\tau^*=0$, $\tau\neq0 $ and
\begin{equation*}
R=\frac{\tau}{4}g\owedge \Bigl( g-4\eta\otimes\eta\Bigr).
\end{equation*}

The latter conclusions support \corref{cor-Fi3-rho} and
\corref{cor-Fi3-R}.

Taking into account \eqref{Rrhotauk-Ex}, we obtain that the
sectional curvatures of the basic $\xi$-sections (respectively,
the basic $\f$-holomorphic sections) are negative (respectively,
 positive)  for $(L,\f, \xi, \eta, g)\in \F_9\oplus\F_{10}$,
which support \thmref{thm-char} for $\F_9$ and $\F_{10}$.



\subsection*{Acknowledgment}
This paper is partially supported by project NI13-FMI-002
of the Scientific Research Fund, Plovdiv University, Bulgaria

\end{document}